# A FUNCTION IN THE NUMBER THEORY

Florentin Smarandache
The University of New Mexico
Department of Mathematics
Gallup, NM 87301, USA

Abstract:
In this paper I shall construct a function1[1] $\eta$ having the following properties:
(1) $\forall\, n \,\varepsilon\, Z,\ n \neq 0,\ (\eta(n))! = M\, n$ (multiple of n).

(2) $\eta(n)$ is the smallest natural number satisfying property (1).

Introduction:
We consider:

$N = \{ 0 , 1 , 2 , 3 , \ldots \}$ and $N^* = \{1, 2, 3, \ldots \}$.

Lemma 1. $\forall\, k, p \,\varepsilon\, N^*,\ p \neq 1$, k is uniquely written

in the form: $k = t_1 a_{n(1)}^{(p)} + \ldots + t_l a_{n(l)}^{(p)}$ where

$$a_{n(i)}^{(p)} = \frac{p^{n(i)} - 1}{p - 1},\ i = \overline{1, l},\ \ n_1 > n_2 > \ldots n_l > 0 \text{ and } 1 \leq t_j \leq p - 1,\ j = \overline{1, l - 1},\ \ 1 \leq t_l \leq p,\ n_i, t_i \,\varepsilon\, N,$$

$i = \overline{1, l}$ , $l \,\varepsilon\, N^*$.

Proof.
The string $( a_n^{(p)})_{n \varepsilon N}$ consists of strictly increasing infinite natural numbers and

$a_{n+1}^{(p)} - 1 = p * a_n^{(p)},\ \alpha n \,\varepsilon\, N^*$, p is fixed,

$a_1^{(p)} = 1,\ a_2^{(p)} = 1 + p,\ a_3^{(p)} = 1 + p + p^2, \ldots$ . Therefore:

$$N^* = \bigcup_{n \varepsilon N^*} ([ a_n^{(p)} , a_{n+1}^{(p)}] \cap N^* ) \text{ where } ( a_n^{(p)} , a_{n+1}^{(p)}) \cap (a_{n+1}^{(p)}, a_{n+2}^{(p)}) = 0$$

because $a_n^{(p)} < a_{n+1}^{(p)} < a_{n+2}^{(p)}$ .

Let $k \,\varepsilon\, N^*$ , $N^* = \bigcup (( a_n^{(p)}, a_{n+1}^{(p)} ) \cap N^* )$,

therefore $\exists!\ n_1 \,\varepsilon\, N^* : k \,\varepsilon\, ( a_{n(1)}^{(p)}, a_{n(1)+1}^{(p)} )$ , therefore k is uniquely written under the form

$$k = \left[ \frac{k}{a_{n_1}^{(p)}} \right] a_{n(1)}^{(p)} + r_1 \text{ (integer division theorem).}$$

[1] This function has been called the Smarandache function. Over one hundred articles, notes, problems and a dozen of books have been written about it.

We note

$$k = \left( \frac{k}{a^{(p)}_{n_1}} \right) = t_1 \rightarrow k = t_1\, a_{n(1)}^{(p)} + r_1,\ r_1 < a_{n(1)}^{(p)}\ .$$

If $r_1 = 0$, as $a_{n(1)}^{(p)} \leq k \leq a_{n(1)+1}^{(p)} - 1 \ \rightarrow\ 1 \leq t_1 \leq p$ and Lemma 1 is proved.

If $r_1 \neq 0$, then $\exists\ !\ n_2\ \varepsilon\ N^* : r_1 [\varepsilon\ \ a_{n(2)}^{(p)}, a_{n(2)+1}^{(p)}]$ ;

$a_{n(1)}^{(p)} > r_1$ involves $n_1 > n_2$, $r_1 \neq 0$ and $a_{n(1)}^{(p)} \leq k \leq a_{n(1)+1}^{(p)} - 1$ involves $1 \leq t_1 \leq p - 1$ because we have

$t_1 \leq (\ a_{n(1)+1}^{(p)}\ - 1 - r_1\ ) : a_n^{(p)}\ < p_1$ .

The procedure continues similarly. After a finite number of steps $l$, we achieve $r_l = 0$, as k = finite, $k\ \varepsilon\ N^*$ and $k > r_1 > r_2 > \ldots > r_l = 0$ and between 0 and k there is only a finite number of distinct natural numbers.

Thus:

k is uniquely written: $k = t_1\, a_{n(1)}^{(p)} + r_1$ , $1 \leq t_1 \leq p - 1$,

r is uniquely written: $r_1 = t_2 * a_{n(2)}^{(p)} + r_2$, $n_2 < n_1$,

$1 \leq t_2 \leq p\text{-}1$,

$r_{l\text{-}1}$ is uniquely written: $r_{l\text{-}1} = t_l * a_{n(l)}^{(p)} + r_l$, and $r_l = 0$,

$n_l < n_{l\text{-}1}$ , $1 \leq t_l \leq p$,

thus k is uniquely written under the form

$k = t_1\, a_{n(1)}^{(p)} + \ldots + t_l\, a_{n(l)}^{(p)}$

with $n_1 > n_2 > \ldots > n_l > 0$, because $n_l\ \varepsilon\ N^*$, $1 \leq t_j \leq p\text{-}1$, $j = \overline{1, l - 1}$, $1 \leq t_l \leq p$, $l \geq 1$.

Let $k\ \varepsilon\ N^*$, $k = t_1 a_{n(1)}^{(p)} + \ldots + t_l a_{n(l)}^{(p)}$ with

$$a_{n(i)}^{(p)} = \frac{p^{ni} - 1}{p - 1} ,$$

$i = \overline{1,\ l}$, $l \geq 1$, $n_i$, $t_i\ \varepsilon\ N^*$, $i = \overline{1, l}$, $n_1 > n_2 > \ldots > n_l > 0$

$1 \leq t_j \leq p - 1$, $j = \overline{1, l - 1}$ , $1 \leq t_l \leq p$.

I construct the function $\eta_p$, p = prime > 0, $\eta_p$: $N^* \rightarrow N$ thus:

$\forall\ n\ \varepsilon\ N^*\ \eta_p(a_n^{(p)}\ ) = p^n$ ,

$\eta_p( t_1 a_{n(1)}^{(p)} + \ldots + t_l a_{n(l)}^{(p)} ) = t_1\ \eta_p(a_{n(1)}^{(p)}) + \ldots + t_l\ \eta_p(a_{n(l)}^{(p)})$.

NOTE 1. The function $\eta_p$ is well defined for each natural number.
Proof

LEMMA 2. $\forall$ k $\varepsilon$ $N^*$, k is uniquely written as $k = t_1 a_{n1}^{(p)} + \ldots + t_l a_{nl}^{(p)}$ with the conditions from Lemma 1, thus $\exists!\ t_1 p^{n(1)} + \ldots + t_l\, p^{n(l)} = \eta_p\,(t_1 a_{n(1)}^{(p)} + \ldots + t_l\, a_{n(l)}^{(p)})$ and $t_1 p^{n(1)} + \ldots + t_l\, p^{n(l)}\ \varepsilon\ N^*$.

LEMMA 3. $\forall$ k $\varepsilon$ $N^*$, $\forall$ p $\varepsilon$ N, p = prime then $k = t_1 a_{n(1)}^{(p)} + \ldots + t_l\, a_{n(l)}^{(p)}$ with the conditions from Lemma 2 thus $\eta_p(k) = t_1 p^{n(1)} + \ldots + t_l\, p^{n(l)}$

It is known that

$$\left[\frac{a_1 + \ldots + a_n}{b}\right] \geq \left[\frac{a_1}{b}\right] + \ldots + \left[\frac{a_n}{b}\right]$$ $\forall$ $a_i$ , b $\varepsilon$ $N^*$ where through [α ] we

have written the integer side of the number α. I shall prove that p's powers sum from the natural numbers which make up the result factors

$( t_1 p^{n(1)} + \ldots + t_l\, p^{n(l)} )\,!$ is $\geq k$;

$$\left[\frac{t_1 p^{n(1)} + \ldots + t_l\, p^{n(l)}}{p}\right] \geq \left[\frac{t_1 p^{n(1)}}{p}\right] + \ldots + \left[\frac{t_l\, p^{n(l)}}{p}\right] =$$

$$t_1 p^{n(1)-1} + \ldots + t_l\, p^{n(l)-1}$$

$$\left[\frac{t_1 p^{n(1)} + \ldots + t_l\, p^{n(l)}}{p^n}\right] \geq \left[\frac{t_1 p^{n(1)}}{p^{n(l)}}\right] + \ldots + \left[\frac{t_l\, p^{n(l)}}{p^{n(l)}}\right] =$$

$$t_1 p^{n(1)-n(l)} + \ldots + t_l\, p^0$$

$$\left[\frac{t_1 p^{n(1)} + \ldots + t_l\, p^{n(l)}}{p^{n(1)}}\right] \geq \left[\frac{t_1 p^{n(1)}}{p^{n(1)}}\right] + \ldots + \left[\frac{t_l\, p^{n(l)}}{p^{n(1)}}\right] =$$

$$t_1 p^0 + \ldots + \frac{t_l\, p^{n(l)}}{p^{n(1)}}\ .$$

Adding → p's powers the sum is $\geq t_1(p^{n(1)-1} + \ldots + p^0) + \ldots + t_l(p^{n(l)-1} + \ldots + p^0)$ =

$t_1\ a_{n(1)}^{(p)} + \ldots + t_l a_{n(l)}^{(p)} = k$.

Theorem 1. The function $n_p$, p = prime, defined previously, has the following properties:

(1) $\exists$ k $\varepsilon$ $N^*$, $(n_p(k))! = M p^k$.

(2) $\eta_p(k)$ is the smallest number with the property (1).

<u>Proof</u>

(1) Results from Lemma 3.

(2) $\forall$ k $\varepsilon$ $N^*$, $p \geq 2$ one has $k = t_1 a_{n(1)}^{(p)} + \ldots + t_l a_{n(l)}^{(p)}$

(by Lemma 2) is uniquely written, where:

$n_i, t_i \, \varepsilon \, N^*$, $n_1 > n_2 > \ldots n_l > 0$,

$$a_{n(i)}^{(p)} = \frac{p^{n(i)} - 1}{p - 1} \; \varepsilon \; N^*,$$

$i = \overline{1, l}$, $1 \leq t_j \leq p - 1$, $j = \overline{1, l-1}$, $1 < t_l < p$.

$\rightarrow \eta_p(k) = t_1 p^{n(1)} + \ldots + t_l p^{n(l)}$ . I note: $z = t_1 p^{n(1)} + \ldots + t_l p^{n(l)}$.

Let us prove that z is the smallest natural number with the property (1). I suppose by the method of reductio ad absurdum that $\exists$ $\gamma$ $\varepsilon$ N, $\gamma < z$:

$\gamma! = M p^k$;

$\gamma < z \rightarrow \gamma \leq z - 1 \rightarrow (z-1)! = M p^k$.

$z - 1 = z = t_1 p^{n(1)} + \ldots + t_l p^{n(l)} - 1$; $n_1 > n_2 > \ldots n_l \geq 1$ and

$n_j \, \varepsilon \, N$, $j = \overline{1, l}$ ;

$$\left[ \frac{z-1}{p} \right] = t_1 p^{n(1)-1} + \ldots + t_{l-1}{}^{n(l-1)-1} + t_l p^{n(l)-1} - 1 \text{ as } \left[ \frac{-1}{p} \right] = -1 \text{ because } p \geq 2,$$

$$\left[ \frac{z-1}{p^{n(l)}} \right] = t_1 p^{n(1) - n(l)} + \ldots + t_{l-1} p^{n(l-1) - n(l)} + t_l p^0 - 1 \text{ as } \left[ \frac{-1}{p^{n(l)}} \right] = -1$$

as $p \geq 2$, $n_l \geq 1$,

$$\left[ \frac{z-1}{p^{n(l)+1}} \right] = t_1 p^{n(1) - n(l) - 1} + \ldots + t_{l-1} p^{n(l-1) - n(l) - 1} + \left[ \frac{t_l p^{n(l)} - 1}{p^{n(l)+1}} \right] =$$

$t_1 p^{n(1) - n(l) - 1} + \ldots + t_{l-1} p^{n(l-1) - n(l) - 1}$ because

$0 < t_l p^{n(l)} - 1 \leq p * p^{n(l)} - 1 < p^{n(l)+1}$ as $t_l < p$;

$$\left[ \frac{z-1}{p^{n(l-1)}} \right] = t_1 p^{n(1) - n(l-1)} + \ldots + t_{l-1} p^0 + \left[ \frac{t_l p^{n(l)} - 1}{p^{n(l-1)}} \right] =$$

$t_1p^{n(1)-n(l-1)} + \ldots + t_{l-1}p^0$ as $n_{l-1} > n_l$ ,

$$\left[\frac{z-1}{p^{n(1)}}\right] = t_1p^0 + \left[\frac{t_2p^{n(2)} + \ldots + t_lp^{n(l)} - 1}{p^{n(1)}}\right] = t_1p^0 .$$

Because $0 < t_2p^{n(2)} + \ldots + t_lp^{n(l)} - 1 \le (p-1)p^{n(2)} + \ldots + (p-1)p^{n(l-1)} + p^*p^{n(l)} - 1 \le$

$$(p-1) * \sum_{i=n(l-1)}^{n_2} p_i + p^{n(l)+1} - 1 \le$$

$$(p-1)\frac{p^{n(2)+1}}{p-1} = p^{n(2)+1} - 1 < p^{n(1)} - 1 < p^{n(1)} \text{ therefore}$$

$$\left[\frac{t_2p^{n(2)} + \ldots + t_lp^{n(l)} - 1}{p^{n(1)}}\right] = 0$$

$$\left[\frac{z-1}{p^{n(1)+1}}\right] = \left[\frac{t_1p^{n(1)} + \ldots + t_lp^{n(l)} - 1}{p^{n(1)+1}}\right] = 0 \text{ because:}$$

$0 < t_1p^{n(1)} + \ldots + t_lp^{n(l)} - 1 < p^{n(1)+1} - 1 < p^{n(1)+1}$ according to a reasoning similar to the previous one.

Adding one gets p's powers sum in the natural numbers which make up the product factors $(z-1)!$ is:

$t_1(p^{n(1)-1} + \ldots + p^0) + \ldots + t_{l-1}(p^{n(l-1)-1} + \ldots + p^0) + t_l(p^{n(l)-1} + \ldots + p^0)$ whence

$1*n_l = k$ or $n_l < k$ or $1 < k$ because

$n_l > 1$ one has $(z-1)! \neq M\,p^k$, this contradicts the supposition made.

Whence $\eta_p(k)$ is the smallest natural number with the property $(\eta_p(k))! = M\,p^k$.

I construct a new function $\eta$: $Z\backslash\{0\} \to N$ defined as follows:

$$\begin{cases} \eta(\pm 1) = 0. \\ \alpha\ n = \varepsilon\, p_1^{\alpha(1)} \ldots p_s^{\alpha(s)} \text{ with } \varepsilon = \pm 1,\ p_i \text{ prime,} \\ p_i = p_j \text{ for } i \neq j,\ \alpha_i \ge 1,\ i = 1, s,\ \eta(n) = \max\limits_{i=1,\ldots,s} \{ \eta_{p_i}(\alpha_i) \}. \end{cases}$$

Note 2. $\eta$ is well defined all over.

Proof

(a) $\forall$ n $\varepsilon$ Z, n $\neq$ 0, n $\neq \pm 1$, n is uniquely written, abstraction of the order of the factors, under the form:

$n = \varepsilon\, p_1^{\alpha(1)} \ldots p_s^{\alpha(s)}$ with $\varepsilon = \pm 1$, where $p_i$ = prime, $p_i \neq p_j$, $\alpha_i \geq 1$ (decomposed into prime factors in Z, which is a factorial ring).

Then $\exists!\ \eta(n) = \max_{i=1,s} \{ \eta_{p(i)}(\alpha_i) \}$ as s = finite and $\eta_{p(i)}(\alpha_i)\ \varepsilon\ N^*$

and $\exists \max_{i=1,\ldots,s} \{\eta_{p(i)}(\alpha_i) \}$

(b) n = $\pm$ 1 $\rightarrow$ E! $\eta(n) = 0$.

Theorem 2. The function $\eta$ previously defined has the following properties:

(1) $(\eta(n))! = M\, n$, $\forall$ n $\varepsilon$ Z\{0};

(2) $\eta(n)$ is the smallest natural number with this property.

Proof

(a) $\eta(n) = \max_{i=1,\ldots,s} \{ \eta_{p(i)}(\alpha_i) \}$, $n = \varepsilon * p_1^{\alpha(1)} \ldots p_s^{\alpha(s)}$ $(n \neq \pm 1)$,

$(\eta_{p(1)}(\alpha_1))! = M\, p_1^{\alpha(1)}$,

$(\eta_{p(s)}(\alpha_s))! = M\, p_s^{\alpha(s)}$.

Supposing $\max_{i=1,\ldots,s} \{ \eta_{p(i)}(a_1) \} = \eta_{p_{i_0}}(\alpha_{i(0)}) \rightarrow (\eta_{p_{i_0}}(\alpha_{i(0)}))! =$

$M\, p_{i(0)}^{\alpha_{i(0)}}$, $\eta_{p_{i_0}}(\alpha_{i_0})\ \varepsilon\ N^*$ and because $(p_i, p_j) = 1$, $i \neq j$,

then $(\eta_{p_{i_0}}(\alpha_{i_0}))! = M\, p_j^{\alpha(j)}$, $\overline{j = 1,s}$ .

Also $(\eta_{p_{i_0}}(\alpha_{i_0}))! = M\, p_1^{\alpha(1)} \ldots p_s^{\alpha(s)}$ .

(b) n = $\pm$ 1 $\rightarrow \eta(n) = 0$; $0! = 1$, $1 = M\, \varepsilon * 1 = M\, n$.

(2) (a) $n \neq \pm 1 \rightarrow n = p_1^{\alpha(1)} \ldots p_s^{\alpha(s)}$ hence $\eta(n) = \max_{i=1,2} \eta_{p(i)}$

Let $\max_{i=1,s} \{ \eta_{p(i)}(\alpha_i) \} = \eta_{p_{i_0}}(\alpha_{i_0})$, $1 \leq i \leq s$;

$\eta_{p_{i_0}}(\alpha_{i_0})$ is the smallest natural number with the property:

$(\eta_{p_{i_0}}(\alpha_{i_0}))! = M\, p_{i_0}^{\alpha_{i(0)}} \rightarrow \alpha\gamma \in N, \gamma < \eta_{p_{i_0}}(\alpha_{i_0})$ whencw

$\gamma! \neq M\, p_{i_0}^{\alpha_{i0}}$ then $\gamma! \neq M \varepsilon * p_1 \ldots p_i^{\alpha_i} \ldots p_s^{\alpha_s} {}_0^{\alpha_{i0}} = M n$ whence

$\eta_{p_{i0}}(\alpha_{i_0})$ is the smallest natural number with the property.

(b) $n = \pm 1 \rightarrow \eta(n) = 0$ and it is the smallest natural number $\rightarrow$ 0 is the smallest natural number with the property $0! = M(\pm 1)$.

NOTE 3. The functions $\eta_p$ are increasing, not injective, on $N^* \rightarrow \{ p^k \mid k = 1, 2, 3, \ldots \}$ they are surjective.

The function $\eta$ is increasing, it is not injective, it is surjective on $Z \setminus \{0\} \rightarrow N \setminus \{1\}$.

CONSEQUENCE. Let $n \in N^*$, $n > 4$. Then n = prime involves $\eta(n) = n$.

Proof

"$\rightarrow$"

n = prime and $n \geq 5$ then $\eta(n) = \eta_n(1) = n$.

"$\leftarrow$"

Let $\eta(n) = n$ and assume by reduction ad absurdum that $n \neq$ prime. Then

(a) $n = p_1^{\alpha(1)} \ldots p_s^{\alpha(s)}$ with $s \geq 2$, $\alpha_i \in N^*$, $i = \overline{1,s}$,

$\eta(n) = \max_{i=1,s} \{ \eta_{p(i)}(\alpha_i) \} = \eta_{p_{i_0}}(\alpha_{i_0}) < \alpha_{i_0} p_{i_0} < n$

contradicting the assumption.

(b) $n = p_1^{\alpha(1)}$ with $\alpha_1 \geq 2$ involves $\eta(n) = \eta_{p(1)}(\alpha_1) \leq p_1 * \alpha_1 < p_1^{\alpha(1)} = n$

because $\alpha_1 \geq 2$ and $n > 4$, which contradicts the hypothesis.

Application

1. Find the smallest natural number with the property:

$n! = M(\pm 2^{31} * 3^{27} * 7^{13})$.

Solution

$\eta(\pm 2^{31} * 3^{27} * 7^{13}) = \max \{ \eta_2(31), \eta_3(27), \eta_7(13) \}$.

Let us calculate $\eta_2(31)$; we make the string

$(a_n^{(2)})_{n \in N^*} = 1, 3, 7, 15, 31, 63, \ldots$

$31 = 1*31 \rightarrow \eta_2(1*31) = 1 * 2^5 = 32$.

Let's calculate $\eta_3(27)$ by making the string

$(a_n^{(3)})_{n\varepsilon N^*}$ = 1, 4, 13, 40, . . . ; 27 = 2*13 + 1 involves $\eta_3(27) = \eta_3(2*13+1*1)$ =

$2*\eta_3(13) + 1*\eta_3(1) = 2*3^3 + 1*3^1 = 54 + 3 = 57$.

Let's calculate $\eta_7(13)$; making the string

$(a_n^{(7)})_{n\varepsilon N^*}$ = 1, 8, 57, . . . ; 13 = 1*8 + 5*1 → $\eta_7(13) = 1*\eta_7(8) + 5*\eta_7(1)$ =

$1*7^2 + 5*7^1 = 49 + 35 = 84 \rightarrow \eta(\pm 2^{31} * 3^{27} * 7^{13}) = \max\{32, 57, 84\} = 84$ involves 84! =

$M(\pm 2^{31} * 3^{27} * 7^{13})$ and 84 is the smallest number with this property.

2. What are the numbers n where n! ends with 1000 zeros?

Solution:

$n = 10^{1000}$, $(\eta(n))! = M\ 10^{1000}$ and it is the smallest number with this property.

$\eta(10^{1000}) = \eta(2^{1000}*5^{1000}) = \max\{\eta_2(1000), \eta_5(1000)\} = \eta_5(1000)$ =

$\eta_5(1*781 + 1*156 + 2*31 + 1) = 1*5^5 + 1*5^4 + 2*5^3 + 1*5^7 = 4005$, 4005 is the smallest

number with this property. 4006, 4007, 4008, 4009 also satisfy this property, but 4010 does not because 4010! = 4009! * 4010 which has 1001 zeros.

Florentin Smarandache 17.11.1979
University of Craiova
Natural Science Faculty
Romania